\documentclass[letterpaper, 10 pt, conference]{ieeeconf}

\IEEEoverridecommandlockouts        
\usepackage{cite}
\usepackage{amsmath,amssymb,amsfonts}
\usepackage{color}
\usepackage{textcomp}
\usepackage{graphicx}
\usepackage{mathrsfs}
\usepackage[vlined,ruled]{algorithm2e}
\usepackage{url}
\usepackage{dsfont}
\usepackage{bbm}
\usepackage{booktabs}
\usepackage{array}
\usepackage[table]{xcolor}
\usepackage{yfonts}
\usepackage{tikz}
\usepackage{multirow}
\usepackage{textcomp}
\usepackage{tabularx}
\usepackage{placeins}
\usepackage{float}
\usepackage{nccmath}
\usepackage{accents}
\usepackage{multicol}
\usepackage{etoolbox, refcount}
\usepackage{chngcntr}
\usepackage{apptools}
\usepackage{relsize}
\usepackage[geometry]{ifsym}
\usepackage[font=footnotesize]{caption}
\usepackage[font=footnotesize]{subcaption}
\usepackage[hidelinks]{hyperref}
\usepackage{cleveref} 
\usepackage{xparse}

\hypersetup{
    colorlinks=true,
    urlcolor=blue
}

\NewDocumentCommand{\INTERVALINNARDS}{ m m }{
    #1 {,} #2
}
\NewDocumentCommand{\interval}{ s m >{\SplitArgument{1}{,}}m m o }{
    \IfBooleanTF{#1}{
        \left#2 \INTERVALINNARDS #3 \right#4
    }{
        \IfValueTF{#5}{
            #5{#2} \INTERVALINNARDS #3 #5{#4}
        }{
            #2 \INTERVALINNARDS #3 #4
        }
    }
}

\newtheorem{theorem}{Theorem}[section]
\newtheorem{lemma}[theorem]{Lemma}
\newtheorem{definition}[theorem]{Definition}
\newtheorem{corollary}[theorem]{Corollary}
\newtheorem{proposition}[theorem]{Proposition}
\newtheorem{problem}{Problem}
\newtheorem{remark}[theorem]{Remark}

\newtheorem{example}[theorem]{Example}

\newcommand{\setdef}[2]{\{#1 \; : \; #2\}}

\newcommand{\Fc}{\mathcal{F}}

\newcommand{\V}{\mathcal{V}}

\newcommand{\real}{\mathbb{R}}

\newcommand{\Cc}{\mathcal{C}}
\newcommand{\Dc}{\mathcal{D}}
\newcommand{\Sc}{\mathcal{S}}
\newcommand{\Pc}{\mathcal{P}}
\newcommand{\Bc}{\mathcal{B}}
\newcommand{\Gc}{\mathcal{G}}
\newcommand{\Vc}{\mathcal{V}}

\newcommand{\Ic}{\mathcal{I}}
\newcommand{\Kc}{\mathcal{K}}

\newcommand{\Nc}{\mathcal{N}}
\newcommand{\Ec}{\mathcal{E}}

\DeclareSymbolFont{bbold}{U}{bbold}{m}{n}
\DeclareSymbolFontAlphabet{\mathbbold}{bbold}

\newcommand{\norm}[1]{\lVert#1\rVert}

\newcommand\oprocendsymbol{\hbox{$\bullet$}}
\newcommand\oprocend{\relax\ifmmode\else\unskip\hfill\fi\oprocendsymbol}

\renewcommand{\theenumi}{(\roman{enumi})}

\newcommand*{\QEDA}{\hfill\ensuremath{\blacksquare}}%

\newcommand\xqed[1]{%
  \leavevmode\unskip\penalty9999 \hbox{}\nobreak\hfill
  \quad\hbox{#1}}
\newcommand\demo{\xqed{$\bullet$}}

\newcommand\problemfinal{\xqed{$\triangle$}}

\newcounter{countitems}
\newcounter{nextitemizecount}
\newcommand{\setupcountitems}{%
  \stepcounter{nextitemizecount}%
  \setcounter{countitems}{0}%
  \preto\item{\stepcounter{countitems}}%
}
\makeatletter
\newcommand{\computecountitems}{%
  \edef\@currentlabel{\number\c@countitems}%
  \label{countitems@\number\numexpr\value{nextitemizecount}-1\relax}%
}
\newcommand{\nextitemizecount}{%
  \getrefnumber{countitems@\number\c@nextitemizecount}%
}
\newcommand{\previtemizecount}{%
  \getrefnumber{countitems@\number\numexpr\value{nextitemizecount}-1\relax}%
}
\makeatother    
{\end{itemize}%
\unskip\computecountitems\ifnumcomp{\previtemizecount}{>}{4}{\end{multicols}}{}}

\newcommand{\longthmtitle}[1]{\mbox{}\emph{(#1):}}

\DeclareFontFamily{U}{mathx}{\hyphenchar\font45}
\DeclareFontShape{U}{mathx}{m}{n}{
<-6> mathx5 <6-7> mathx6 <7-8> matha7
<8-9> mathx8 <9-10> mathx9
<10-12> mathx10 <12-> mathx12
}{}
\DeclareSymbolFont{mathx}{U}{mathx}{m}{n}
\DeclareMathSymbol{\bigplus}{\mathop}{mathx}{"90}
\DeclareMathSymbol{\bigtimes}{\mathop}{mathx}{"91}

\renewcommand{\baselinestretch}{.96}
\newcommand{\comment}[1]{} 
\newcolumntype{P}[1]{>{\centering\arraybackslash}p{#1}}
\allowdisplaybreaks

\begin{document}
\title{\bf Stabilization of Nonlinear Systems through Control Barrier Functions}

\author{Pol Mestres \qquad Kehan Long \qquad Melvin Leok \qquad Nikolay Atanasov \qquad Jorge Cort{\'e}s%
\thanks{P. Mestres, K. Long, M. Leok, N. Atanasov and J. Cort{\'e}s are with the Contextual Robotics Institute, University of California, San Diego (e-mails: \{pomestre,k3long,mleok,natanasov,cortes\}@ucsd.edu).
}
}

\maketitle

\begin{abstract}
This paper proposes a control design approach for stabilizing nonlinear control systems. 
Our key observation is that the set of points where the decrease condition of a control Lyapunov function (CLF) is feasible can be regarded as a safe set.
By leveraging a nonsmooth version of control barrier functions (CBFs) and a weaker notion of CLF, we develop a control design that forces the system to converge to and remain in the region where the CLF decrease condition is feasible.
We characterize the conditions under which our controller asymptotically stabilizes the origin or a small neighborhood around it, even in the cases where it is discontinuous. We illustrate our design in various examples.
\end{abstract}

\section{Introduction}

Control Lyapunov functions (CLFs)~\cite{EDS:98} are a well-established tool for designing stabilizing controllers for nonlinear systems. CLF-based control designs ensure that the controller satisfies a Lyapunov decrease condition, guaranteeing asymptotic stability of the origin. However, finding a CLF for a general nonlinear control system is challenging, even though sum-of-squares~\cite{WT:06} or neural network~\cite{YCC-NR-SG:19} techniques have been proposed. On the other hand, control barrier functions (CBFs)~\cite{PW-FA:07,ADA-SC-ME-GN-KS-PT:19} are widely used in safety-critical applications to design controllers that enforce predefined safety specifications.
Boolean nonsmooth control barrier functions (BNCBFs) \cite{PG-JC-ME:17-csl,PG-IB-ME:19} extend CBF theory to a richer class of safe sets that cannot be expressed as the superlevel set of a differentiable function. In the setting where both safety and stability must be certified, several works have proposed approaches to combine CLFs and CBFs \cite{MZR-BJ:16, KL-CQ-JC-NA:21-ral, ADA-SC-ME-GN-KS-PT:19,PM-JC:23-csl}.
The key novel idea that we explore in this paper is that it is often possible to construct candidate CLFs for which the Lyapunov decrease condition is feasible in large regions of the state space even if they are not valid CLFs. Instead of modifying this candidate CLF to be a valid CLF, we consider the set of points where the Lyapunov decrease condition is feasible as a \emph{safe set}. This brings up the question of whether the notion of CLF can be relaxed and combined with CBFs to yield a design methodology for stabilizing controllers. 
The ideas in this paper are related to \cite{YC-MJ-MS-ADA:21}, which extends the safe operating region of a controller by implementing a \textit{backup controller}, and 
\cite{ART-NK:97}, which devises a control strategy that combines local and global stabilizing controllers.


\subsubsection*{\textbf{Statement of Contributions}}
This work considers stabilizing nonlinear control systems. First, we introduce the notion of Weak Control Lyapunov Function (WCLF), which relaxes the Lyapunov decrease condition to be feasible only in a subset of the state space that need not include an open neighborhood of the origin.
Next, we interpret the set where the Lyapunov decrease condition is feasible as a \emph{safe set} and use BNCBFs to design a controller that keeps the system within the safe set while satisfying the Lyapunov decrease condition. If the BNCBF condition is feasible outside the safe set, we extend our control strategy to ensure trajectories starting outside achieve finite-time convergence to the safe set.
Our result shows that Filippov solutions of the closed-loop system (coinciding with standard solutions if the controller is continuous) with an initial condition in the safe set asymptotically converge to the smallest sublevel set of the WCLF that does not contain \emph{incompatible} points outside it, i.e., points where the Lyapunov decrease condition and the BNCBF condition can not be satisfied simultaneously.
Lastly, we showcase our control design's applicability in three examples.
For reasons of space, proofs are omitted and will appear
elsewhere.

\section{Preliminaries}

We introduce preliminaries on discontinuous dynamical systems, weak control Lyapunov functions, and Boolean nonsmooth control barrier functions.

\subsubsection*{Notation}
We denote by $\mathbb{Z}_{>0}$, $\real$ and $\real_{\geq0}$ the set of positive integers, real numbers, and non-negative real numbers, respectively. Given a set $\Sc\subset\real^n$, we write $\text{Int}(\Sc)$, $\partial \Sc$, $\overline{co}(\Sc)$ for the interior, the boundary and the convex closure of $\Sc$, respectively. The $n$-dimensional zero vector is denoted by $\textbf{0}_n$, and $\norm{x}$ denotes the Euclidean norm of $x\in\real^n$. 
For $\delta>0$ and $x\in\real^n$, we let $\Bc(x,\delta)=\setdef{y\in\real^n}{\norm{y-x}\leq \delta}$.
Given $f:\real^n\to\real^n$, $g:\real^n\to\real^{n\times m}$ and a smooth function $W:\real^n\to\real$, the Lie derivatives of $W$ with respect to $f$ and $g$ are $L_fW = \nabla W^T f$ and $L_gW = \nabla W^T g$, respectively. 
A function $\beta:\real\to\real$ is of extended class $\Kc$ if $\beta(0)=0$ and $\beta$ is strictly increasing. A function $V:\real^n\to\real$ is positive-definite if $V(0)=0$ and $\forall x\neq0, \; V(x)>0$.
Let $F:\real^n\to\real^n$ be a locally Lipschitz vector field and consider the system $\dot{x} = F(x)$. Local Lipschitzness of $F$ ensures that, for every initial condition $x_0\in\real^n$, there exists $T>0$ and a unique trajectory $x(t;x_0)$ such that $\frac{d}{dt}x(t;x_0) = F(x(t;x_0))$ for all $t\in[0,T]$ and $x(0;x_0)=x_0$. 
A set $\Pc$ is forward-invariant if $x_0\in\Pc$ implies $\forall t\geq0, \; x(t;x_0)\in\Pc$. If $\Pc$ is forward-invariant and $x^*\in\Pc$ is an equilibrium, $x^*$ is Lyapunov stable relative to $\Pc$ if for every open set $U$ containing $x^*$, there exists an open set $\tilde{U}$ also containing $x^*$ such that for all $ x_0\in\tilde{U}\cap\Pc$,
%
%
$x(t;x_0)\in U\cap\Pc$ for all $t>0$. The equilibrium $x^*$ is asymptotically stable relative to $\Pc$ if it is Lyapunov stable relative to $\Pc$ and there is an open set $U$ containing $x^*$ such that $\lim\limits_{t\to\infty} x(t;x_0) = x^*$ for all $ x_0 \in U\cap\Pc$.
Given a locally Lipschitz function $h:\real^n\to\real$, the generalized gradient of $h$ at $x\in\real^n$ is $\partial h(x) \!=\! \overline{co}\setdef{ \lim\limits_{i\to\infty} \nabla h(x_i) }{ x_i\to x, x_i\notin S\cup\Omega_f }$, where $\Omega_f$ is the zero-measure set where $f$ is nondifferentiable and $S$ can be any set of measure zero. 

\subsubsection*{Discontinuous Dynamical Systems}
Consider the differential equation 
\begin{align}\label{eq:discontinuous-ode}
\dot{x} = F(x),
\end{align}
where $F:\real^n\to\real^n$ is measurable and essentially locally bounded (cf.\cite{AFF:88}). For $x\in\real^n$, let
$\Fc[F](x) = \bigcap\limits_{\delta>0} \bigcap\limits_{\mu(S) = 0} \overline{co} \{ F(\Bc(x,\delta)\backslash S) \}$, where $\mu(S)$ is the Lebesgue measure of $S$.
A Filippov solution of~\eqref{eq:discontinuous-ode} on $[t_0,t_1]\subset\real$ is a solution of the differential inclusion $\dot{x} \in \Fc[F](x)$, i.e., an absolutely continuous function $[t_0,t_1]\to\real^n$ such that $\dot{x}(t)\in \Fc[F](x)$ for almost all $t\in[t_0,t_1]$.

\subsubsection*{Weak Control Lyapunov Functions and Strict Boolean Nonsmooth Control Barrier Functions}
Consider a control-affine system
\begin{align}\label{eq:control-affine-sys}
  \dot{x}=f(x)+g(x)u,
\end{align}
where $f:\real^{n}\to\real^{n}$ and $g:\real^{n}\to\real^{n\times m}$
are locally Lipschitz functions, with $x\in\real^{n}$ the state and
$u\in\real^{m}$ the input. Throughout the paper, and without loss of generality, we assume $f(\textbf{0}_n)=\textbf{0}_n$, so that the origin $\textbf{0}_n$ is the desired equilibrium state of the (unforced) system. 

\begin{definition}\longthmtitle{Weak Control Lyapunov Function}\label{def:cclf}
  Given an open set $\mathcal{D}\subseteq\real^{n}$, with $\textbf{0}_n \in \Dc$, a continuously differentiable function
  $V:\real^{n}\to\real$ is a weak control Lyapunov function (WCLF) in $\mathcal{D}$ for the
  system~\eqref{eq:control-affine-sys} if
  $V$ is proper in $\mathcal{D}$, i.e.,
  $\setdef{x\in\mathcal{D}}{V(x)\leq c}$ is a compact set for all $c>0$, $V$ is positive-definite, and there exists a continuous positive-definite function
  $W:\real^{n}\to\real$ and a set $\tilde{\Dc}\subset\Dc$ such that, for each $x\in\mathcal{\tilde{\Dc}}$, there exists a control $u\in\real^{m}$ satisfying
  \begin{align}\label{eq:clf-ineq}
    L_fV(x)+L_gV(x)u \leq -W(x).
  \end{align}
\end{definition}
\medskip

If $\tilde{\Dc}$ in Definition~\ref{def:cclf} is an open set containing the origin, then the notion of WCLF is equivalent to CLF~\cite{EDS:98,RAF-PVK:96a}. 
If $V$ is a CLF, any Lipschitz controller $\hat{u}:\real^n\to\real^m$ that satisfies~\eqref{eq:clf-ineq} for all
$x\in\mathcal{D}$ asymptotically stabilizes the origin~\cite{EDS:98}.
However, the set $\tilde{\Dc}$ in Definition~\ref{def:cclf} need not include the origin. WCLFs guarantee the existence of a control law that decreases the value of $V$ for all points in $\tilde{\Dc}$, but such control law does not guarantee the asymptotic stabilization of the origin because it might steer the system towards states outside of $\tilde{\Dc}$.

Next, we define the notion of strict Boolean nonsmooth control barrier function (SBNCBF), adapted from~\cite{PG-JC-ME:18-ccta,PG-JC-ME:17-csl}.

\begin{definition}\longthmtitle{Strict Boolean Nonsmooth Control Barrier Function}\label{def:strict-bncbf}
    Let $N\in\mathbb{Z}_{>0}$, and let $h_i:\real^n\to\real$, $i\in[N]$, be continuously differentiable functions.
    Let $h(x)=\max_{i\in[N]}h_i(x)$ and
    \begin{align}\label{eq:set-c}
      \Cc = \setdef{x \! \in \! \real^n \! \!}{\! h(x) \! \geq \! 0}, \; \;
      \partial \Cc = \setdef{x \! \in \! \real^n \! \! }{\! h(x) \! = \! 0}.
    \end{align}
    We also let the set of active constraints at $x$ be $\Ic(x):=\setdef{i\in[N]}{h(x)=h_i(x)}$.
    The function $h:\real^n\to\real$
    is a strict Boolean nonsmooth control barrier function (SBNCBF) of $\Cc$ if there exists an open set $\Gc\subset\real^n$ containing $\Cc$, an extended class $\Kc$ function $\alpha:\real\to\real$ and $\epsilon>0$ such that for all $x\in\Gc$
    there exists a neighborhood $\Nc_x$ of $x$ such that for all $y\in\Nc_x$
    there exists $u\in\real^m$ satisfying
    \begin{align}\label{eq:bncbf-ineq}
        \min_{v\in\partial h(x)} v^T (f(y)+g(y)u) \geq -\alpha(h(y)) + \epsilon.
    \end{align}
\end{definition}

When $N=1$ and $\epsilon=0$, Definition~\ref{def:strict-bncbf} reduces to the standard notion of CBF~\cite[Definition 2]{ADA-SC-ME-GN-KS-PT:19}. 
However, SBNCBFs allow for a richer class of safe sets, which motivates their use in this work.
Moreover, if $h$ is a SBNCBF,~\cite[Theorem 3]{PG-JC-ME:17-csl} shows that if there exists a Lipschitz controller $\hat{u}:\real^n\to\real^m$ and a neighborhood $\Nc_x$ of every $x\in\Cc$ such that $\hat{u}$ satisfies~\eqref{eq:bncbf-ineq} for all $y\in\Nc_x$ then $\hat{u}$ makes $\Cc$ forward invariant.
The requirement that~\eqref{eq:bncbf-ineq} is satisfied with $\epsilon>0$ is necessary for some of the results in the paper.

The following result, adapted from~\cite[Theorem 3]{PG-JC-ME:18-ccta}, provides a sufficient condition for $h$ to satisfy Definition~\ref{def:strict-bncbf}.
\begin{proposition}\longthmtitle{Sufficient condition for SBNCBF}\label{prop:suff-cond-bncbf}
Suppose there exists an open set $\Gc\subset\real^n$ containing $\Cc$ and a locally Lipschitz extended class $\Kc$ function $\alpha:\real\to\real$ and $\epsilon>0$ such that for all $x\in\Gc$ there exists a neighborhood $\Nc_x$ of $x$, and for all $y\in\Nc_x$ there exists $u\in\real^m$ that satisfies
  \begin{align}\label{eq:bncbf-suff-cond}
    L_f h_i(y) + L_g h_i(y)u \geq -\alpha(h_i(y)) + \epsilon,
  \end{align}
  for all $i\in\Ic(x)$. Then, $h$ is a SBNCBF of $\Cc$.
\end{proposition}
\smallskip

When dealing with both safety and stability specifications, we note 
that an input $u$ might satisfy~\eqref{eq:clf-ineq} but not~\eqref{eq:bncbf-suff-cond}, or vice versa. 
The following notion captures when both constraints can be satisfied simultaneously and is adapted from~\cite{PM-JC:23-csl}.
\begin{definition}\longthmtitle{Compatibility of WCLF-SBNCBF pair}\label{def:clf-bncbf-compatibility}
  Let $\mathcal{D} \subseteq \real^{n}$ be open, $\Cc \subset \Dc$ be
  closed, $V$ a WCLF in $\mathcal{D}$ and $h$ a SBNCBF of~$\Cc$. Then, $V$
  and $h$ are a compatible WCLF-SBNCBF pair at $x \in \Cc$ if there exists a neighborhood $\Nc_x$ of $x$ 
  such that for all $y\in\Nc_x$ there exists
  $u\in\real^{m}$ satisfying~\eqref{eq:clf-ineq} at $y$
  and~\eqref{eq:bncbf-suff-cond}
  for all $i\in\Ic(x)$ simultaneously. We refer to both functions
  as a compatible WCLF-SBNCBF pair in a set $\Ec \subset \real^n$ if $V$ and $h$ are a compatible WCLF-SBNCBF pair at every point 
  in $\Ec$.
\end{definition}
\smallskip

\section{Problem Statement}\label{sec:pb-statement}

Consider a control-affine system of the form~\eqref{eq:control-affine-sys}.
Let $V:\real^n\to\real$ be a WCLF for~\eqref{eq:control-affine-sys} on a set $\Dc\subset\real^n$ and suppose that the set $\tilde{\Dc} \subset \Dc$ in Definition~\ref{def:cclf} is known. We consider the following problem.


\begin{problem}\label{pb:problem-1}
  Find a control law $\bar{u}:\real^n\to\real^m$
  and a region $\Gamma\subset\real^n$ such that
  trajectories of~\eqref{eq:control-affine-sys} with initial condition in $\Gamma$ asymptotically converge to the origin.
  \problemfinal
\end{problem}

The key insight to solve this problem is that the set $\tilde{\Dc}$ can be treated as a \textit{safe set} (because~\eqref{eq:clf-ineq} is feasible at $\tilde{\Dc}$). Hence, if we can find a set $\Cc\subset\tilde{\Dc}$, and a SBNCBF of $\Cc$ that is compatible with $V$ in $\Cc$, then we can define a control law $\bar{u}_1$ that steers the system trajectories towards the origin and remains in $\Cc$. Moreover, if the SBNCBF is feasible outside of $\Cc$, we can extend $\bar{u}_1$ (potentially discontinuously) 
so that it steers trajectories outside of $\Cc$ towards it. Hence the set $\Cc \subset\tilde{\Dc}\subset\Dc$ can be used to construct $\Gamma$ in Problem~\ref{pb:problem-1}.

\section{Stabilizing Control Design using WCLFs and SBNCBFs}\label{sec:proposed-solution}

\begin{figure}[htb]
  \centering
  %
  \includegraphics[width=0.45\textwidth]{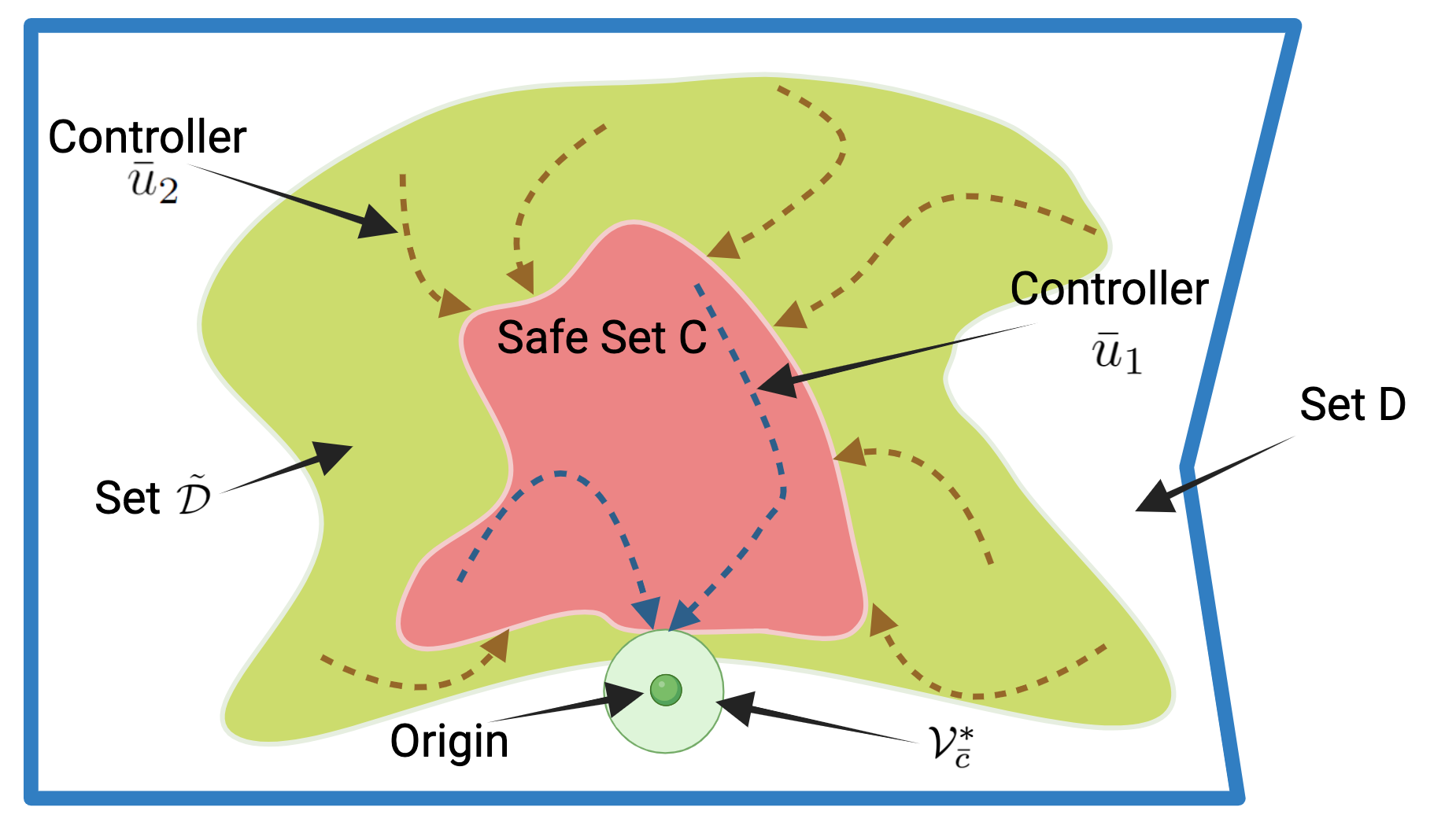}
  \caption{Illustration of the control design in Proposition~\ref{prop:invariance-convergence-practical}.}\label{fig:control_illustrate}
  \vspace*{-3ex}
\end{figure}

This section formalizes our control design idea to solve Problem~\ref{pb:problem-1}. Let $h_i:\real^n\to\real$, $i\in[N]$, $N\in\mathbb{Z}_{>0}$, be continuously differentiable functions, $h(x)=\max_{i\in[N]}h_i(x)$ and define $\Cc$ as in~\eqref{eq:set-c}.
Suppose that 
$\Cc$ is connected and $\Cc\subset\tilde{\Dc}$.
For any $c>0$, let $\Vc_c^*:=\setdef{x\in\real^n}{V(x)\leq c}$.
Next we present the main result of the paper, which solves Problem~\ref{pb:problem-1} and is illustrated in Figure~\ref{fig:control_illustrate}.

%
%

\begin{proposition}\longthmtitle{Invariance and convergence to smallest compatible Lyapunov level set}\label{prop:invariance-convergence-practical}
    Suppose that $h$ is a SBNCBF of $\Cc$.
    Let $\bar{c}>0$ be such that $\Vc_{\bar{c}}^*\cap\Cc\neq\emptyset$ and suppose that $V$ and $h$ are a compatible WCLF-SBNCBF pair in $\Cc\backslash\Vc_{\bar{c}}^*$.
    Let $\bar{u}_1:\real^n\to\real^m$ be a locally Lipschitz controller such that
    \begin{enumerate}
        \item $\bar{u}_1$ satisfies~\eqref{eq:clf-ineq} for all $x\in\Cc\backslash\Vc_{\bar{c}}^*$;
        \item for all $x\in\Cc\backslash \Vc_{\bar{c}}^*$
        there exists a neighborhood $\Nc_x$ of $x$ such that $\bar{u}_1$ satisfies~\eqref{eq:bncbf-suff-cond} for all $y\in\Nc_x$ and $i\in\Ic(x)$.
    \end{enumerate}
    Moreover, let $\Gc\subset\real^n$ 
     be the set where the SBNCBF is feasible (cf. Definition~\ref{def:strict-bncbf}), and suppose that $\Gc$ is connected.
    Let $\bar{u}_2:\real^n\to\real^m$ be a locally Lipschitz controller such that for all $x\in(\Gc\backslash\Cc) \cup(\Cc\cap\Vc_{\bar{c}}^*)$, there exists a neighborhood $\Nc_x$ of $x$ such that $\bar{u}_2$ satisfies~\eqref{eq:bncbf-suff-cond} for all $y\in\Nc_x$ and $i\in\Ic(x)$.
    Define
    \begin{align*}
        \bar{u}(x) = \begin{cases}
            \bar{u}_1(x) \quad \text{if} \quad x\in\Cc\backslash\Vc_{\bar{c}}^*, \\
            \bar{u}_2(x) \quad \text{otherwise},
        \end{cases}
    \end{align*}
    and consider the closed-loop system
    \begin{align}\label{eq:closed-loop-discontinuous-practical-stab}
        \dot{x}=f(x)+g(x)\bar{u}(x).
    \end{align}
    Then,~\eqref{eq:closed-loop-discontinuous-practical-stab} has a unique Filippov solution $\bar{x}(t;x_0)$ from any initial condition $x_0$. Moreover, for any $c>\bar{c}$,
    \begin{enumerate}
        \item\label{it:third} if 
        $x_0\in\Vc_c^*\cap\Cc$, then $\bar{x}(t;x_0)\in\Vc_c^*\cap\Cc$ for all $t\geq0$;
        \item\label{it:fourth} if $x_0\in\Cc\backslash\Vc_c^*$, then there exists $t_{1}>0$ such that $\bar{x}(t_{1};x_0)\in\Vc_c^*\cap\Cc$ and $\bar{x}(t;x_0)\in\Vc_c^*\cap\Cc$ for all $t\geq t_1$;
        \item\label{it:fifth} if $x_0\notin\Cc$ and $\bar{x}(t;x_0)\in\Gc$ for all $t\geq0$, then there exists $t_2>0$ such that $\bar{x}(t_2;x_0)\in\Cc$ and $\bar{x}(t;x_0)\in\Cc$ for all $t\geq t_2$. Moreover, there exists $t_3\geq t_2$ such that $\bar{x}(t_3;x_0)\in\Vc_c^*\cap\Cc$ and $\bar{x}(t;x_0)\in\Vc_c^*\cap\Cc$ for all $t\geq t_3$.
    \end{enumerate}
\end{proposition}
\smallskip

The following result specializes Proposition~\ref{prop:invariance-convergence-practical} to the case where the origin is in $\Cc$ and $V$ and $h$ are a compatible WCLF-SBNCBF pair in $\Cc$.

\begin{corollary}\longthmtitle{Invariance and convergence to the origin}\label{cor:invariance-convergence-origin}
    Suppose that $h$ is a SBNCBF of $\Cc$.
    Further suppose that $V$ and $h$ are a compatible WCLF-SBNCBF pair in $\Cc$ and $\textbf{0}_n \in \Cc$. Let $\Gc\subset\real^n$ be as in Definition~\ref{def:strict-bncbf}.
    Take $\bar{c}=0$ and define $\bar{u}_1:\real^n\to\real^m$, $\bar{u}_2:\real^n\to\real^m$ and $\bar{u}:\real^n\to\real^m$ as in Proposition~\ref{prop:invariance-convergence-practical}.
    Then,~\eqref{eq:closed-loop-discontinuous-practical-stab} has a unique Filippov solution $\bar{x}(t;x_0)$ from any initial condition $x_0\in\real^n$. Moreover,
    \begin{enumerate}
        \item\label{it:first} if $x_0\in\Cc$, then $\bar{x}(t;x_0)\in\Cc$ for all $t\geq0$ and $\lim\limits_{t\to\infty}\bar{x}(t;x_0)=\textbf{0}_n$;
        \item\label{it:second} if $x_0\notin\Cc$ and $\bar{x}(t;x_0)\in\Gc$ for all $t\geq0$, then there exists $t_4>0$ such that $\bar{x}(t_4;x_0)\in\Cc$ and $\bar{x}(t;x_0)\in\Cc$ for all $t\geq t_4$. Moreover, $\lim\limits_{t\to\infty}\bar{x}(t;x_0) = \textbf{0}_n$.
    \end{enumerate}
\end{corollary}
\smallskip

Leveraging Proposition~\ref{prop:invariance-convergence-practical} and Corollary~\ref{cor:invariance-convergence-origin}, our control design methodology takes the following steps.
\begin{enumerate}
    \item[{(1)}] Find a WCLF and identify the set $\tilde{\Dc}$; 
    \item[{(2)}] Find a set $\Cc \subset \tilde{\Dc} \subset \real^n$  and a SBNCBF $h$ of $\Cc$; 
    \item[{(3)}] Find a sublevel set $\Vc_{\bar{c}}^*$ of $V$ such that $V$ and $h$ are a compatible WCLF-SBNCBF pair in $\Cc\backslash\Vc_{\bar{c}}^*$.
\end{enumerate}

\begin{remark}\longthmtitle{Classical solutions}\label{rem:classical-sols}
    In general, the controller $\bar{u}$ in Proposition~\ref{prop:invariance-convergence-practical} and Corollary~\ref{cor:invariance-convergence-origin} is discontinuous. However, if
    $\bar{u}$ is locally Lipschitz, then the results in Proposition~\ref{prop:invariance-convergence-practical} and Corollary~\ref{cor:invariance-convergence-origin} hold with classical (instead of Filippov) solutions. Indeed, if $f+g\bar{u}$ is continuous at $x\in\real^n$ then the set $\Fc[f+g\bar{u}](x)$ is equal to $f(x)+g(x)\bar{u}(x)$ and Filippov solutions coincide with classical ones.
    \demo
\end{remark}

\begin{remark}\longthmtitle{Conditions on $\Gc$}\label{rem:h-not-global-bncbf}
    Proposition~\ref{prop:invariance-convergence-practical} requires that $x(t;x_0)\in\Gc$ for all $t\geq0$. In general, this condition is difficult to verify. However, this condition holds
    if $\Gc=\real^n$ or $\Gc$ is a superlevel set of $h$.
    \demo
\end{remark}

\begin{remark}\longthmtitle{Stability of the origin}\label{rem:lyap-stab-origin}
    Under the assumptions in Corollary~\ref{cor:invariance-convergence-origin}:
    \begin{enumerate}
        \item if the origin is in $\text{Int}(\Cc)$, Corollary~\ref{cor:invariance-convergence-origin} guarantees that the origin is asymptotically stable;
        \item if the origin is in $\partial\Cc$, Corollary~\ref{cor:invariance-convergence-origin} guarantees that the origin is asymptotically stable \textit{relative to} $\Cc$. However, in this case Lyapunov stability of the origin is not guaranteed, since trajectories that start close to $\Cc$ but outside of it might take a long excursion away from the origin before
        entering $\Cc$ and converging to it;

        \item the origin can not be outside of $\Cc$. Indeed, Proposition~\ref{prop:invariance-convergence-practical} guarantees that we can design a controller that makes all trajectories with initial condition in $\Cc$ stay in $\Cc$ for all future times 
        and always decrease the value of $V$, which is not possible if the origin is not in $\Cc$.
        \demo
    \end{enumerate}
\end{remark}

\begin{remark}\longthmtitle{Lipschitz controller with relaxed CLF condition}\label{rem:lipschitz-relaxed-clf}
    If the controller $\bar{u}$ in Proposition~\ref{prop:invariance-convergence-practical} and Corollary~\ref{cor:invariance-convergence-origin} cannot be designed continuously, we give the following alternative design.
    Let $\breve{u}:\real^n\to\real^m$ be a locally Lipschitz controller satisfying~\eqref{eq:bncbf-suff-cond} and the relaxed version of~\eqref{eq:clf-ineq}
        $L_fV(x) + L_gV(x) \breve{u}(x) + W(x) \leq \delta(x)$,
    where $\delta:\real^n\to\real_{\geq0}$. 
    For example, given $\lambda>0$, one can take 
    \begin{align}\label{eq:relaxed-cclf-bncbf-qp}
        \breve{u}(x) = &\min_{u\in\real^m, \delta\in\real} \frac{1}{2}\norm{u}^2 + \lambda \delta^2, \\
        \notag
        &\quad \text{s.t.} \quad \quad ~\eqref{eq:bncbf-suff-cond}, \ L_fV(x) + L_gV(x) u + W(x) \leq \delta.
    \end{align}
    The work~\cite{PM-AA-JC:24-ejc} gives conditions under which $\breve{u}$ is locally Lipschitz. 
    Even though $\breve{u}$ has no stability guarantees because the CLF condition is relaxed, in Section~\ref{sec:examples} we show how this controller has good performance properties in practice. \demo
\end{remark}


\begin{remark}\longthmtitle{Compact safe sets}\label{rem:compact-safe-set}
    Even though Proposition~\ref{prop:invariance-convergence-practical} does not require $\Cc$ to be compact, verifying its assumptions is often easier if $\Cc$ is compact (for example, using Lemma~\ref{lem:boundary-sufficient}).
    If $\Cc$ is not compact, we can take $\gamma\in\real_{\geq0}$ large enough so that
    $\Cc\cap\Vc_{\gamma}^*\neq\emptyset$ and consider a new compact safe set defined by $\tilde{\Cc}_{\gamma} = \Cc\cap\Vc_{\gamma}^*$.
    Moreover, if $V$ is a WCLF in $\Cc$ and $V$ and $h$ are a compatible WCLF-SBNCBF pair in $\Cc$, then it follows that $\tilde{h}(x) = \min\{ h(x), \gamma-V(x) \}$ is an SBNCBF of $\tilde{\Cc}_{\gamma}$ and $V$ and $\tilde{h}$ are a compatible WCLF-SBNCBF pair in $\tilde{\Cc}_{\gamma}$.
    \demo
\end{remark}


\section{Illustrative Examples}\label{sec:examples}

We demonstrate our control design in three examples.\footnote{Open-source implementations of the examples are available at \href{https://github.com/KehanLong/CBF_Stabilization}{https://github.com/KehanLong/CBF\_Stabilization}.}

\begin{figure*}[htb]
\centering
\subfloat[Nominal WCLF QP Controller]
{\includegraphics[width=0.33\textwidth]{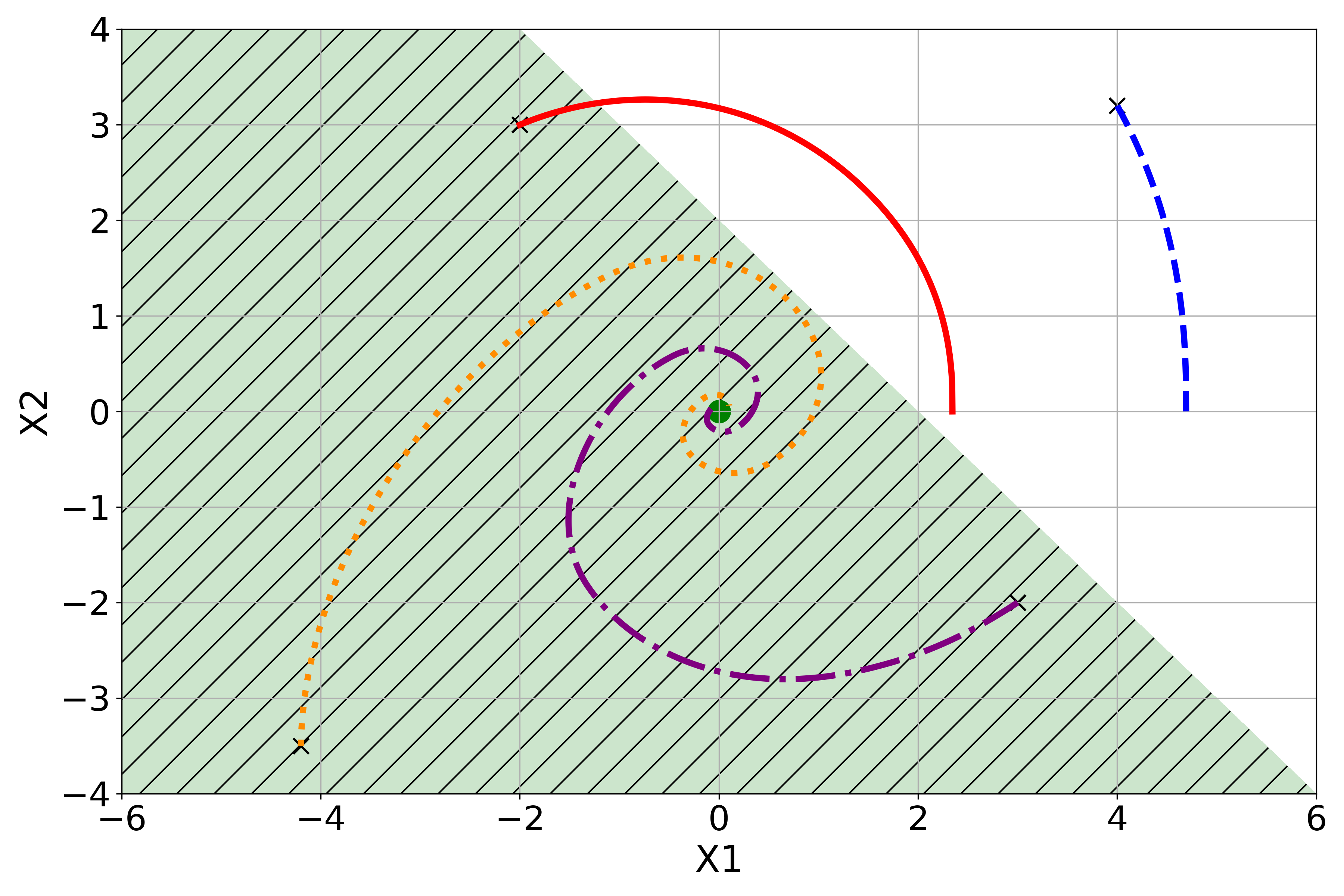}
\label{fig:warm_up_clf_only}}
\hfill
\subfloat[Switching WCLF-SBNCBF QP Controller.]{\includegraphics[width=0.33\textwidth]{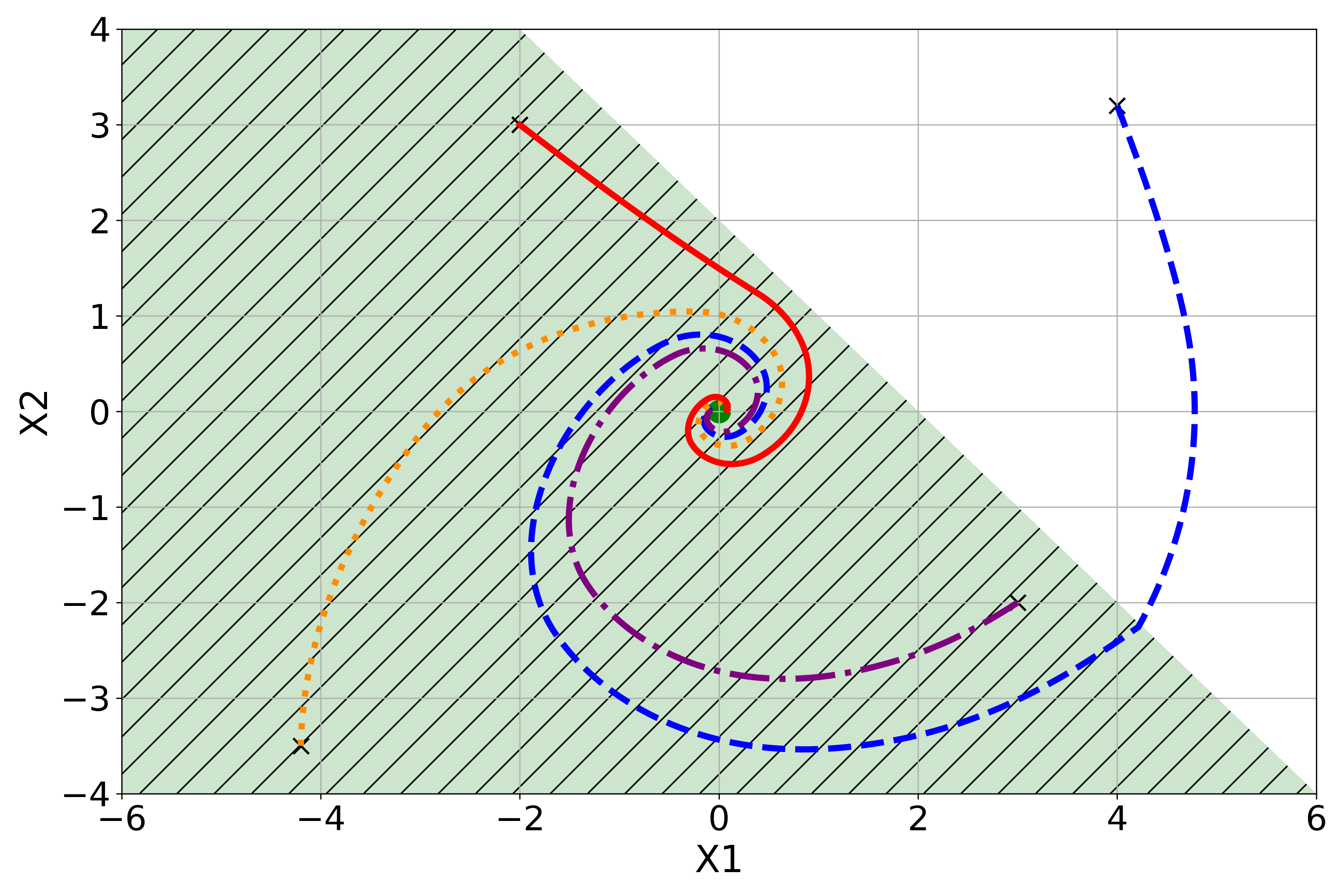}\label{fig:clf_cbf_switch_warmup}}
\hfill
\subfloat[Relaxed WCLF-SBNCBF QP Controller.]{\includegraphics[width=0.33\textwidth]{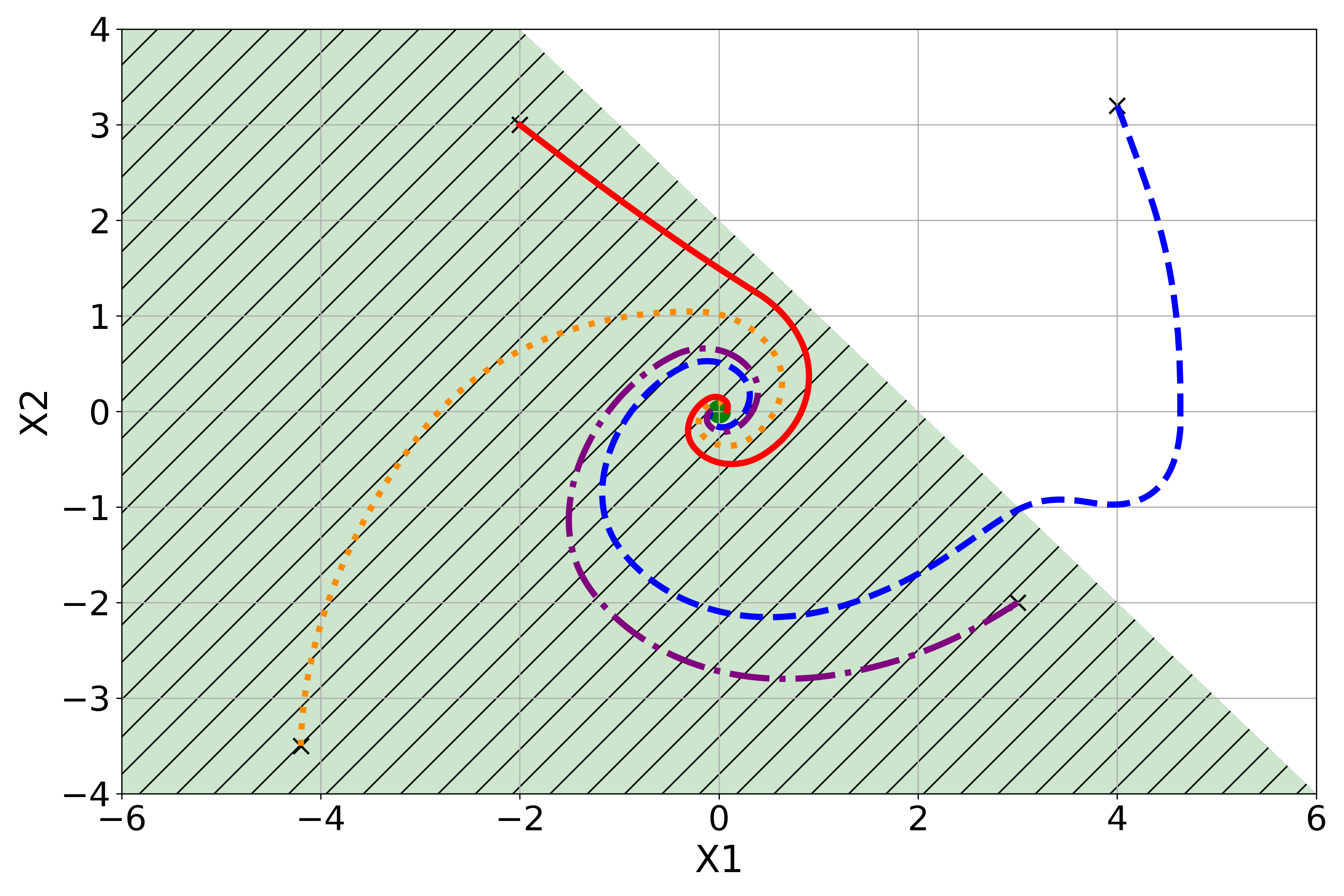}\label{fig:clf_cbf_relax_warmup}}\\
\caption{Comparison of controller performances: WCLF QP, switching WCLF-SBNCBF QP, and relaxed WCLF-SBNCBF QP controllers for the example in Example~\ref{sec:warmup-example}. The safe region defined by the BNCBF is depicted in light green. Each of the four initial states is marked as a black cross and the system's equilibrium is shown as a green dot. In Fig.~\ref{fig:warm_up_clf_only}, only relying on the quadratic WCLF fails to stabilize the system to the equilibrium, as the CLF condition cannot be satisfied once $x_2 = 0$ and $x_1 \geq 2.02$. On the other hand, with the proposed SBNCBF, either the switching WCLF-SBNCBF QP or the relaxed WCLF-SBNCBF QP controller effectively stabilizes the system to equilibrium for all four initial states.}
\label{fig:warmup_trajectories}
\vspace{-3ex}
\end{figure*}

\begin{figure}[htb]
  \centering
  %
  \includegraphics[width=0.4\textwidth]{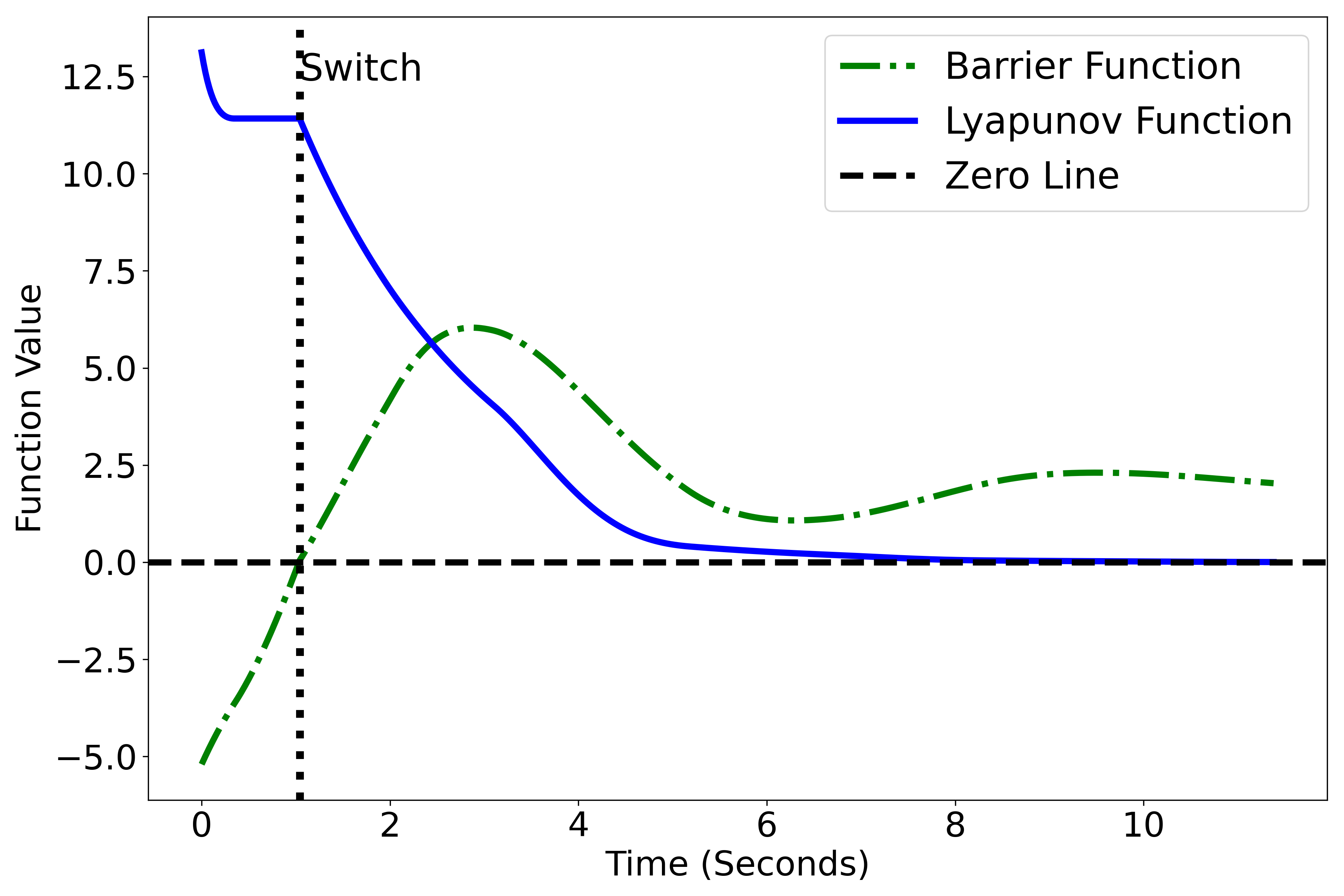}
  \caption{Functions values over time for Trajectory 1 in Fig.~\ref{fig:clf_cbf_switch_warmup}.  This illustrates the values of SBNCBF (green) and WCLF (blue) throughout the trajectory. }\label{fig:warmup_values}
  \vspace*{-3ex}
\end{figure}

\begin{example}\label{sec:warmup-example}
Let $\Lambda:\real\to\real$ be defined as $\Lambda(s) = e^{-\frac{1}{2.02-s}}$ if $s < 2.02$ and $\Lambda(s)=0$ if $s \geq 2.02$.
We note that $\Lambda$ is locally Lipschitz
and consider the dynamics
\begin{align}\label{eq:warmup-dynamics}
    \dot{x}_1 = x_2 - x_1 \Lambda(x_1), \quad \dot{x}_2 = -x_1 + u,
\end{align}
with $x = [x_1, x_2] \in \real$ the state and $u\in\real$ the input. The function $V_w:\real^2\to\real$, $V_w(x) = \frac{1}{2}(x_1^2 + x_2^2)$, is a WCLF. Since $\dot{V}_w(x)\!=\!-x_1^2 \Lambda(x_1) \!+\! x_2 u$, the set $\tilde{\Dc}$ in Definition~\ref{def:cclf} can be taken as $\tilde{\Dc}_w\!=\!\real^2\backslash\setdef{x\in\real^2}{x_2=0, x_1 \geq 2.01}$. Now, define $h_w:\real^2\to\real$ as $h_w(x) = -x_1-x_2+2$ and let $\Cc_w := \setdef{x\in\real^2}{h_w(x)\geq0}\subset\tilde{\Dc}_w$.
Note that $\Cc_w$ is not compact but for any $\gamma>0$ we can define a compact subset of it as
$\tilde{\Cc}_{w,\gamma} = \setdef{x\in\real^2}{V_w(x)\leq \gamma} \cap \Cc_w$, using the construction in Remark~\ref{rem:compact-safe-set}.
Note that $\tilde{\Cc}_{w,\gamma} = \setdef{x\in\real^2}{\tilde{h}_{w,\gamma}(x)=\min\{ h_w(x), \gamma-V_w(x) \} \geq 0 }$.
Next, we show that for any $\gamma>0$, $\tilde{h}_{w,\gamma}$ is a SBNCBF of $\tilde{\Cc}_{w,\gamma}$ and $(V_w,\tilde{h}_{w,\gamma})$ is a compatible WCLF-SBNCBF pair in $\tilde{\Cc}_{w,\gamma}$.

\underline{\textbf{$\tilde{h}_{w,\gamma}$ is a SBNCBF of $\tilde{\Cc}_{w,\gamma}$}}.
Since $\Cc_{w}$ is only defined by a single continuously differentiable function, and $\dot{h}_w(x_1,x_2)=-x_2+x_1 \Lambda(x_1) + x_1 - u$, 
for all $(x_1,x_2)\in\partial\tilde{\Cc}_{w,\gamma}$ with $h_w(x_1,x_2)=0$ and $V_w(x)\neq \gamma$, there exists a neighborhood $\Nc_x$ of $x$ such that~\eqref{eq:bncbf-suff-cond} is feasible for all points in $\Nc_x$ for any $\epsilon>0$ and $\alpha$.
At points $x\in\partial\tilde{\Cc}_{w,\gamma}$ where $V(x)=\gamma$, and $h_w(x)\neq0$,~\eqref{eq:bncbf-suff-cond} is feasible in a neighborhood $\Nc_x$ of $x$ because $V_w$ is a WCLF. Finally, if $x$ is such that $V(x)=\gamma$ and $h_w(x)=0$, the fact that~\eqref{eq:bncbf-suff-cond} is feasible in a neighborhood of $x$ follows from the fact that $V_w$ and $\tilde{h}_{w,\gamma}$ are a compatible WCLF-SBNCBF pair, which we show next.
Since the SBNCBF condition is feasible for all points in $\partial\tilde{\Cc}_{w,\gamma}$, Lemma~\ref{lem:boundary-sufficient} ensures that $\tilde{h}_{w,\gamma}$ is a SBNCBF of $\tilde{\Cc}_{w,\gamma}$.

\underline{\textbf{$V_w$ and $\tilde{h}_{w,\gamma}$ are a compatible WCLF-SBNCBF pair}}
\underline{\textbf{in $\tilde{\Cc}_{w,\gamma}$}}.
Let $x\in\partial\Cc\cap\tilde{\Cc}_{w,\gamma}$. If $x_2 > 0$, there exists a  neighborhood $\Nc_x$ of $x$ and $u\in\real$ sufficiently negative and large in absolute value such that~\eqref{eq:clf-ineq} and~\eqref{eq:bncbf-suff-cond} are simultaneously feasible for all points in $\Nc_x$ for any $\epsilon>0$, extended class $\Kc$ function $\alpha$ and positive definite function $W$.
If $x\in\partial\Cc\cap\tilde{\Cc}_{w,\gamma}$ and $x_2=0$, there exists a sufficiently small neighborhood $\Nc_x$ of $x$ and a positive definite function $W(x)=\sigma_0(x_1^2+x_2^2)$ with $\sigma_0>0$ sufficiently small such that any $u$ that satisfies~\eqref{eq:bncbf-suff-cond} for $y\in\Nc_x$, also satisfies~\eqref{eq:clf-ineq} at $y$. Hence, there exists $\tilde{\delta}$ sufficiently small such that $V_w$ and $h_w$ are a compatible WCLF-SBNCBF pair for all $x\in\partial\Cc\cap\tilde{\Cc}_{w,\gamma}$ with $x_2\in(-\tilde{\delta},0)$.
Finally, take $x\in\partial\Cc\cap\tilde{\Cc}_{w,\gamma}$, with $x_2 \leq -\tilde{\delta}$,
$\epsilon < 2$,
let $M_\gamma=\sup_{x\in\tilde{\Cc}_{w,\gamma}} x_1^2+x_2^2$, and take $W(x) = \sigma(x_1^2+x_2^2)$ with $\sigma < \min\{ \frac{2\tilde{\delta}^2}{M_\gamma}, \sigma_0 \}$. It follows that $-x_2+x_1\Lambda(x_1)+x_1 - \epsilon \geq \frac{x_1^2 \Lambda(x_1)-\sigma( x_1^2 + x_2^2) }{x_2}$, which implies that $V_w$ and $\tilde{h}_{w,\gamma}$ are a compatible WCLF-SBNCBF pair at $x$ (cf.~\cite[Lemma 5.2]{PM-JC:23-csl}).
Hence, we have proved that $V_w$ and $\tilde{h}_{w,\gamma}$ are a compatible WCLF-SBNCBF pair in $\partial\Cc\cap\tilde{\Cc}_{w,\gamma}$. Now, the fact that $V_w$ and $\tilde{h}_{w,\gamma}$ are a compatible WCLF-SBNCBF pair in all of $\partial\tilde{\Cc}_{w,\gamma}$ follows from the fact that $V_w$ is a WCLF and by constructing a linear extended class $\Kc$ function with a similar argument as in Lemma~\ref{lem:boundary-sufficient}.
Now, by using Lemma~\ref{lem:boundary-sufficient} it follows that for any $\gamma>0$, $\tilde{h}_{w,\gamma}$ is a SBNCBF of $\tilde{\Cc}_{w,\gamma}$ and the set $\Gc$ in Definition~\ref{def:strict-bncbf} can be taken as $\Gc=\real^n$. Moreover, $V_w$ and $\tilde{h}_{w,\gamma}$ are a WCLF-SBNCBF compatible pair in $\tilde{\Cc}_{w,\gamma}$. Therefore, the results in Corollary~\ref{cor:invariance-convergence-origin} apply.

\underline{\textbf{Simulation results}.} In the simulation, we set $W$ as $W(x) = 0.5V_w(x)$ and $\alpha$ as $\alpha(s)=s$ for all $s\in\real$. In Fig.~\ref{fig:warmup_trajectories}, we compare the performance of the controller obtained as the solution of the quadratic program (QP) that at every state minimizes the norm of the controller and satisfies only~\eqref{eq:clf-ineq} (denoted as WCLF-QP), the controller obtained as the solution of the QP that at every state minimizes the norm of the controller and satisfies~\eqref{eq:clf-ineq} and~\eqref{eq:bncbf-suff-cond} for all $x\in\Cc_w$ and satisfies the assumptions of Proposition~\ref{prop:invariance-convergence-practical} (denoted as switching WCLF-SBNCBF-QP) 
and the relaxed WCLF-SBNCBF-QP controller (presented in Remark~\ref{rem:lipschitz-relaxed-clf}) for four initial states. In Fig.~\ref{fig:warm_up_clf_only}, the WCLF-QP fails to stabilize the system for two initial states, since~\eqref{eq:clf-ineq} cannot be satisfied once $x_2=0$ and $x_1\geq2.02$. On the other hand, both the switching WCLF-SBNCBF-QP and the relaxed WCLF-SBNCBF-QP controller successfully stabilize the system. 
Moreover, when the system is outside the safe set, the satisfaction of the SBNCBF condition drives the system to the safe set, leading to a temporary non-decrease in Lyapunov function values. At around $t=1$ seconds, the system enters the safe set, and the satisfaction of the WCLF and SBNCBF conditions ensures stabilization to the origin without leaving the safe set.
\end{example}

\begin{figure}[t]
  \centering
  \subcaptionbox{WCLF QP Controller\label{fig:unicycle_clf_only}}{\includegraphics[width=0.46\linewidth]{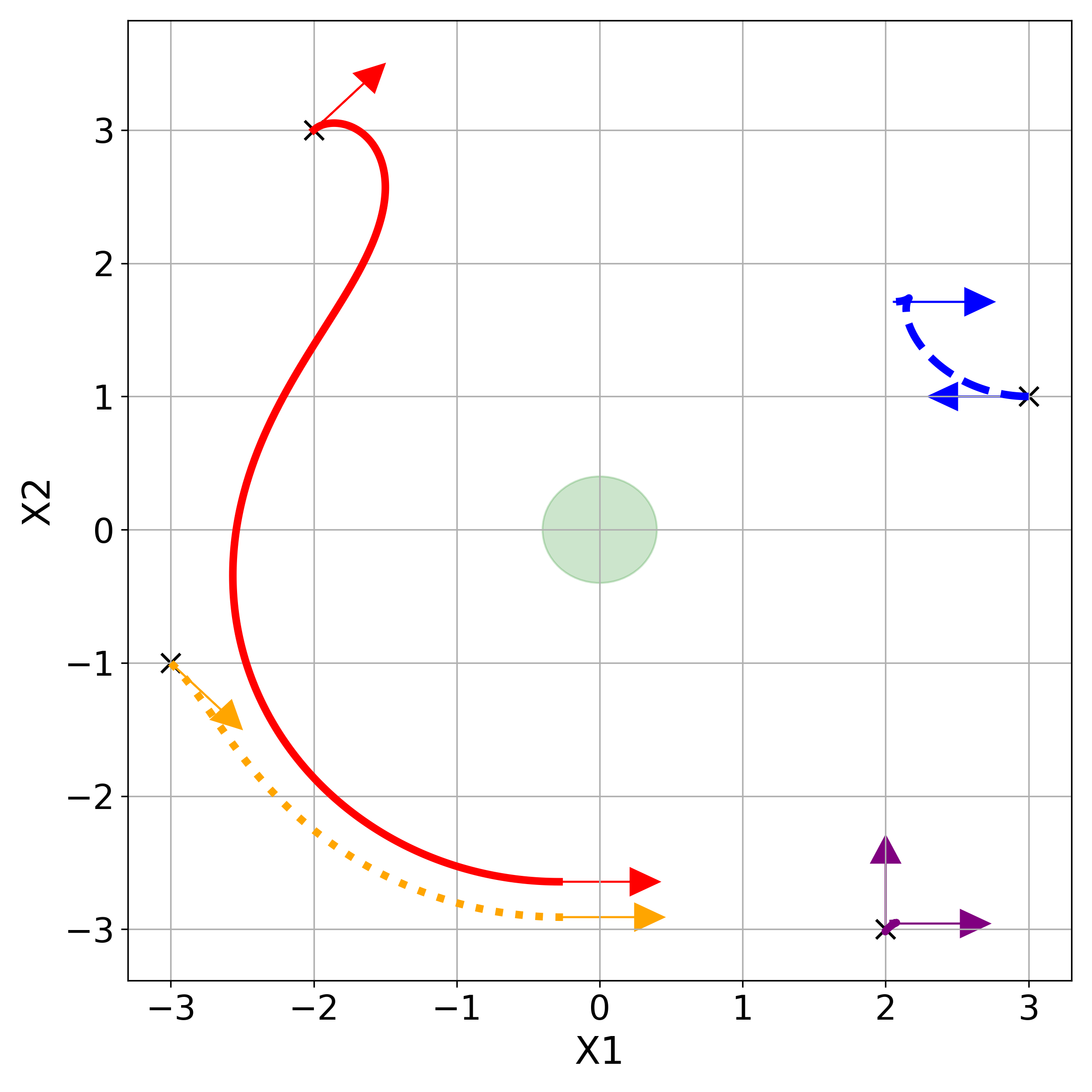}}%
  \hfill%
  \subcaptionbox{Switching WCLF-SBNCBF QP Controller\label{fig:unicycle_switching}}{\includegraphics[width=0.46\linewidth]{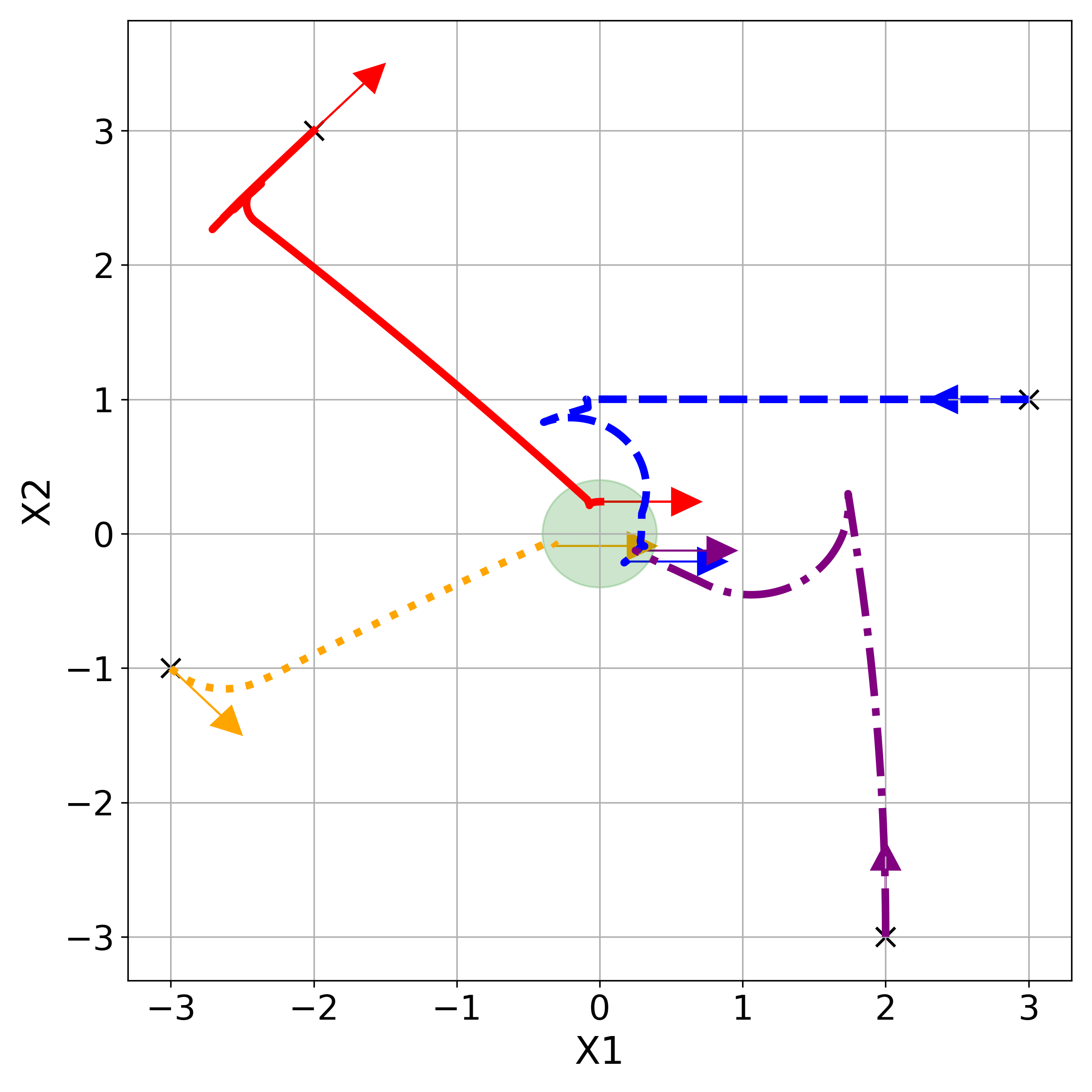}}\\
  \caption{Comparison of controller performances: WCLF QP and switching WCLF-SBNCBF QP for a unicycle system in Example~\ref{sec:unicycle-dynamics-example}. The designated goal region is highlighted in light green, which is a disk centered at $(0,0)$ with a radius of $r=0.4$ and contains the set of points where compatibility cannot be verified. The black crosses represent the initial states and the arrows indicate the system's orientation at the start and end of each trajectory, showcasing the final alignment with the desired orientation. }
  \label{fig:unicycle}
  \vspace*{-3ex}
\end{figure}

\begin{example}\label{sec:unicycle-dynamics-example}

Consider unicycle dynamics:
\begin{equation*}\label{eq:unicycle}    
\dot{x}_1 = v \cos(x_3), \quad \dot{x}_2 = v \sin(x_3), \quad \dot{x}_3 = w,
\end{equation*}
with state $x=[x_1,x_2,x_3]\in\real^3$ and inputs $v, w\in\real$. We consider stabilizing the system at the origin but our approach can be adapted to stabilize the system at any point in $\real^3$. Consider the WCLF $V_u(x)=x_1^2 + x_2^2 + bx_3^2$,
with $b>0$ a parameter to be designed. Let $W_u:\real^3\to\real$ be the associated positive definite function in Definition~\ref{def:cclf}.
Note that $V_u$ is not a CLF and the CLF condition~\eqref{eq:clf-ineq} reads
\begin{align}\label{eq:clf-condition-unicycle}
    &2x_1 v \cos(x_3) + 2x_2 v \sin(x_3) + 2b x_3 w + W(x) \leq 0.
\end{align}
If $x_3=0$ and $x_1\cos(x_3)+x_2\sin(x_3) = 0$, then \eqref{eq:clf-condition-unicycle} cannot be satisfied unless $x_1 = x_2 = 0$. Therefore, the set $\tilde{\Dc}$ in Definition~\ref{def:cclf} can be taken as 
\begin{align*}
    &\tilde{\Dc}_{u} := \real^3\backslash \setdef{x\in\real^3}{x_3 = 0, \ x_1\cos(x_3)+x_2\sin(x_3) = 0}.
\end{align*}
Now, let $\delta>0$ and define
\begin{align}
\label{eq:bncbf-unicycle}
    &h_{1,u}(x) = \delta-( -x_1\sin(x_3) + x_2\cos(x_3) )^2, \notag\\
    &h_{2,u}(x) = x_1^2 + x_2^2 - 1.5^2 \delta, \notag\\
    &h_u(x) = \min \{ h_{1,u}(x), h_{2,u}(x) \}.
\end{align}
Further, let $\Cc_{u}\!:=\!\setdef{x\in\real^3}{h_{u}(x)\geq 0}$,
which is connected, but not compact. 
Let $\tilde{\Cc}_u$ be a compact subset of $\Cc_u$ obtained using the compactification procedure described in Remark~\ref{rem:compact-safe-set} and followed in the previous example. Following an argument similar to the previous example, it is sufficient to show that $h_u$ satisfies the SBNCBF condition at $\partial\Cc_u$ and that $V_u$ and $h_u$ are a compatible WCLF-SBNCBF pair in $\tilde{\Cc}_u\backslash\V_{\bar{c}}^*$, where $\V_{\bar{c}}^*$ is a small sublevel set of $V_u$ to be designed.
We first show that the set $\Cc_u$ (and hence also $\tilde{\Cc}_u$) is a subset of $\tilde{\Dc}_u$, as required in Definition~\ref{def:cclf}.
Throughout this example we let $a(x)=-2(-x_1 \sin(x_3)+x_2\cos(x_3))$.
and $\bar{a}(x) = x_1\cos(x_3)+x_2\sin(x_3)$.

\underline{\textbf{The set inclusion $\Cc_{u} \subset \tilde{\Dc}_{u}$ holds.}}
Suppose that $x\in\Cc_{u}$ and $x\notin\tilde{\Dc}_{u}$.
Then, $x_3 \!=\! 0$ and
$\bar{a}(x)=0$.
Since $(x_1,x_2)$ and $(\cos(x_3), \sin(x_3) )$ are orthogonal and $(-\sin(x_3),\cos(x_3))$ and $(\cos(x_3), \sin(x_3) )$ are also orthogonal, $(x_1,x_2)$ is proportional to $(-\sin(x_3), \cos(x_3) )$. By the Cauchy-Schwartz inequality, this means that $|-x_1\sin(x_3) + x_2\cos(x_3)| = \sqrt{ x_1^2 + x_2^2 }$.
Now, if $h_{1,u}(x) \geq 0$, $\delta \geq x_1^2 + x_2^2$.
This means that $h_{2,u}(x) < 0$ and thus $z\notin\Cc_{u}$, reaching a contradiction.

\underline{\textbf{$h_u$ is a SBNCBF of $\Cc_u$.}}
We show that there exists $\epsilon>0$ and an extended class $\Kc$ function $\alpha$ such that for all $x\in\partial\Cc_u$, there exists a neighborhood $\Nc_x$ such that~\eqref{eq:bncbf-suff-cond} is feasible.
First, suppose that $h_{1,u}(x) = 0$ and $h_{2,u}(x)>0$.
Condition~\eqref{eq:bncbf-suff-cond} at $x$ for $h_{1,u}$ reads
\begin{align}\label{eq:strict-bncbf-unicycle-1}
    &2a(x)\bar{a}(x)\omega + \alpha(h_{1,u}(x)) \geq \epsilon.
\end{align}
Note that $a(x) = \pm 2\sqrt{\delta}\neq0$.
Moreover, $\bar{a}(x) \neq 0$.
Indeed, if $\bar{a}(x) = 0$, then
by the same argument used to show that $\Cc_u\subset\tilde{\Dc}_u$ we have $|-x_1\sin(x_3)+x_2\cos(x_3)| = \sqrt{ x_1^2 + x_2^2 } > 1.5\sqrt{\delta}$,
where in the last inequality we have used $h_{2,u}(x) > 0$. This contradicts $h_{1,u}(x)=0$. Hence, $\bar{a}(x) \neq 0$ and, if $h_{1,u}(x) = 0$ and $h_{2,u}(x) > 0$, there exists a neighborhood $\Nc_x$ of $x$ for which~\eqref{eq:strict-bncbf-unicycle-1} is feasible at all points in $\Nc_x$ for any $\epsilon$ and $\alpha$.
Next, suppose $h_{2,u}(x) = 0$ and $h_{1,u}(x)>0$.
Condition~\eqref{eq:bncbf-suff-cond} at $z$ for $h_{2,u}$ reads
\begin{align}\label{eq:strict-bncbf-unicycle-2}
    2\bar{a}(x) v \! + \! \alpha(h_{2,u}(x)) \! \geq \! \epsilon.
\end{align}
If $\bar{a}(x) = 0$, by the same argument used to show $\Cc_u\subset\tilde{\Dc}_u$, we have that
$|\!-\!x_1\sin(x_3)\!+\!x_2\cos(x_3)| = \sqrt{ x_1^2 + x_2^2 } \!=\! 1.5\sqrt{\delta}$, which contradicts $h_{1,u}(x) \!<\! 0$.
Hence, there exists a neighborhood $\Nc_x$ of $x$ for which~\eqref{eq:strict-bncbf-unicycle-2} is feasible for all points in $\Nc_x$ for any $\epsilon$ and $\alpha$.
Lastly, if $h_{1,u}(x) \!=\! h_{2,u}(x) = 0$, since~\eqref{eq:strict-bncbf-unicycle-1}
can be satisfied using only $\omega$ and~\eqref{eq:strict-bncbf-unicycle-2}
can be satisfied using only $v$, there also exists a neighborhood $\Nc_x$ of $x$ for which~\eqref{eq:strict-bncbf-unicycle-1} and~\eqref{eq:strict-bncbf-unicycle-2} are simultaneously feasible 
for all points in $\Nc_x$ for any $\epsilon$ and $\alpha$.

\underline{\textbf{Compatible region for $V_u$ and $h_u$}.}
Let $\bar{c}=4\delta$ and $B\!>\!0$ such that $|x_3-x_3|<\sqrt{B}$ for all $x\in\tilde{\Cc}_u$ (which exists because $\tilde{\Cc}_u$ is compact), and take $b=\delta/B$.
Using the notation in Section~\ref{sec:proposed-solution}, we show that $V_u$ and $h_u$ are a compatible WCLF-SBNCBF pair in $\tilde{\Cc}_u \backslash \Vc_{\bar{c}}^*$.
First, we show that $V_u$ and $h_{1,u}$ are compatible in $\tilde{\Cc}_u \backslash \Vc_{\bar{c}}^*$.
We use~\cite[Lemma 5.1]{PM-JC:23-csl}, which gives a characterization of when a CLF
and a CBF are compatible at a point. Let $x\in\tilde{\Cc}_u \backslash \Vc_{\bar{c}}^*$ and suppose that there exists $\kappa>0$ such that
\begin{align}
\label{eq:compat-cond-ld}
  2\bar{a}(x) = 0, \quad x_3 = -\kappa a(x)\bar{a}(x).
\end{align}
Then, $x_3=0$, which implies that the WCLF condition for $V_u$ only involves $v$ and the SBNCBF condition for $h_{1,u}$ only involves $\omega$. Therefore, inequalities~\eqref{eq:clf-ineq} for $V_u$ and~\eqref{eq:bncbf-suff-cond} for $h_{1,u}$ can be satisfied simultaneously in a neighborhood of $x$ for any $\epsilon$ and $\alpha$.
  This implies that $V_u$ and $h_{1,u}$ are a compatible WCLF-SBNCBF pair in $\tilde{\Cc}_u\backslash \Vc_{\bar{c}}^*$.
  Now, let $k_u:\real^3\to\real^2$ be a continuous controller satisfying the WCLF condition for $V_u$ and the SBNCBF condition for $h_{u,1}$ for all $x\in\tilde{\Cc}_u\backslash \Vc_{\bar{c}}^*$. Such controller exists by~\cite[Proposition 3.1]{PO-JC:19-cdc}, by taking $\epsilon>0$ sufficiently small.
  Since $b = \delta/B $ and $\bar{c}=4\delta$,
  $h_{2,u}(x)\geq 0.5 \delta$ for all $x\in\tilde{\Cc}_u\backslash \Vc_{\bar{c}}^*$.
  Now, consider the linear extended class $\Kc$ function $\alpha_2:\real\to\real$, $\alpha_2(s) = \alpha_{2,0}s$ with
  \begin{align}\label{eq:alpha-20}
      \alpha_{2,0} > \frac{1}{0.5 \delta }\sup\limits_{x\in\tilde{\Cc}_u\backslash \Vc_{\bar{c}}^*}\left| \begin{pmatrix}
          \bar{a}(x) \\ 
          0
      \end{pmatrix}^T k_u(x) - \epsilon\right|,
  \end{align}
  the right hand side of~\eqref{eq:alpha-20} is bounded because $k_u$ is continuous and $\tilde{\Cc}_u$ is compact. Using $\alpha_2$ as extended class $\Kc$ function, $k_u$ satisfies the WCLF condition for $V_u$ and the SBNCBF condition for $h_{u}$. Hence, $V_u$ and $h_u$ are a compatible WCLF-SBNCBF pair in $\tilde{\Cc}_u$.
  Hence, the assumptions of Proposition~\ref{prop:invariance-convergence-practical} hold and our control design ensures that all trajectories that start in $\tilde{\Cc}_u$ converge to $\tilde{\Cc}_u\cap\Vc_{\bar{c}}^*$.
  The fact that $h_u$ satisfies the SBNCBF condition at $\partial\Cc_u$ ensures that there exists a set $\Gc$ containing $\tilde{\Cc}_u$ as in Definition~\ref{def:strict-bncbf}. 
  However,~\eqref{eq:bncbf-suff-cond} is in general not feasible outside of $\tilde{\Cc}_u$.

\underline{\textbf{Simulation results}.} In the simulation, we specify $\delta=0.04$, $b = 0.01$, positive-definite function $W(x) = 0.1V(x)$, and define the extended class $\mathcal{K}$ function $\alpha(s) = 0.005 s$. 
As shown in Fig.~\ref{fig:unicycle_clf_only}, relying solely on the WCLF fails to stabilize the system to the origin, since trajectories end up at points where the CLF condition is not feasible. However, as shown in Fig.~\ref{fig:unicycle_switching}, the switching WCLF-SBNCBF QP controller 
converges to a neighborhood of the origin $\Vc_{\bar{c}}^*$.
\end{example}

\section{Conclusions}

This paper proposed an approach for stabilizing nonlinear systems using weak CLFs that are not valid CLFs. Our key idea is to treat the subset where the CLF condition is feasible as a safe set and utilize a non-smooth CBF to keep the system trajectories in the safe set. We proved that the proposed controller has stability guarantees both when it is continuous or discontinuous, using appropriate notions of solution for the closed-loop system. 
Our methodology requires the identification of a WCLF-SBNCBF pair, together with a set where both are compatible. We have illustrated this process in different examples.
Future work will focus on three fronts. Firstly, we aim to develop theoretical and computational tools to simplify the process of identifying a compatible WCLF-SBNCBF pair. Secondly, we plan to investigate explicit control designs that satisfy the requirements in our results and ensure continuity of the resulting controller.
Thirdly, we plan to extend our method to systems with uncertainty. 


\section*{Acknowledgements}
The authors are grateful to M. Alyaseen for multiple conversations on barrier functions for discontinuous systems.

\bibliography{bib/alias,bib/JC,bib/Main-add,bib/Main,bib/New}
\bibliographystyle{IEEEtran}

\appendix

\begin{lemma}\longthmtitle{Checking SBNCBF or compatibility conditions on the boundary is sufficient}\label{lem:boundary-sufficient}
    Suppose that $\Cc$ is a compact set.
    \begin{enumerate}
        \item\label{it:boundary-first} Suppose that there exists $\epsilon>0$ and a neighborhood $\Nc_x$ for all $x\in\partial\Cc$, such that for all $y\in\Nc_x$ there exists $u\in\real^m$ and $\epsilon>0$ satisfying
        \begin{align}\label{eq:boundary-bncbf-condition}
            L_fh_i(y)+L_gh_i(y)u \geq \epsilon,
        \end{align}
        for all $i\in\Ic(x)$. Then, $h$ is a SBNCBF of $\Cc$.
        \item\label{it:boundary-second} Suppose that there exists $\epsilon>0$, a positive definite function $W:\real^n\to\real$ and a neighborhood $\Nc_x$ of every $x\in\partial\Cc$ such that 
        for all $y\in\bar{\Nc}_x$ there exists $u\in\real^m$ satisfying~\eqref{eq:boundary-bncbf-condition} for all $i\in\Ic(x)$ and~\eqref{eq:clf-ineq} at $y$. Then, $V$ and $h$ are a compatible WCLF-SBNCBF pair in $\Cc$.
    \end{enumerate}
\end{lemma}

\end{document}